\documentclass[11]{amsart}

\usepackage{amsfonts}
\usepackage{amscd}
\usepackage{amsmath, mathrsfs, amssymb, mathtools}
\usepackage{amsthm}
\usepackage{setspace}
\usepackage{hyperref}
\usepackage{color}
\usepackage{epsfig}
\usepackage{here}
\usepackage{graphicx}
\usepackage[all]{xy}
\usepackage{psfrag}
\usepackage{graphicx,transparent}
\usepackage{enumerate}
\usepackage{caption}

\theoremstyle{plain}
\newtheorem{theorem}{Theorem}[section]

\newtheorem{problem}[theorem]{Problem}
\theoremstyle{definition}

\begin{document}

\title[Topology of moduli spaces of differentials]{An algebro-geometric perspective on the topology of moduli spaces of differentials}

\date{\today}

\author{Dawei Chen}
\address{Department of Mathematics, Boston College, Chestnut Hill, MA 02467, USA}
\email{dawei.chen@bc.edu}

\author{Fei Yu}
\address{School of Mathematical Sciences, Zhejiang University, Hangzhou, 310058, China}
\email{yufei@zju.edu.cn}

\thanks{Research of D.C. is supported in part by the National Science Foundation under Grant DMS-2301030 and by Simons Travel Support for Mathematicians.}

\begin{abstract}
    Differentials on Riemann surfaces correspond to translation surfaces with conical singularities, and affine transformations acting on them preserve the orders of these singularities. This viewpoint allows the moduli spaces of differentials to appear in various guises across many areas, including algebraic geometry, dynamical systems, combinatorial enumeration, and mathematical physics.  Over the past few decades, remarkable progress has been made in computing invariants of these moduli spaces, classifying linear subvarieties, understanding degenerations and compactifications, and developing intersection theory on these spaces.  Despite these advances, our understanding of the topology of moduli spaces of differentials remains limited, and many fundamental questions are still open.  In this survey, we aim to present, from an algebro-geometric perspective, the known results and open problems concerning the topology of moduli spaces of differentials, as well as their connections to other aspects of the field, with the hope of inspiring further developments in the coming decade.
\end{abstract}

\maketitle

\setcounter{tocdepth}{1}
\tableofcontents

\section{Introduction}
\label{sec:intro}

Holomorphic differentials on Riemann surfaces correspond to translation surfaces with conical singularities, where the zero orders of the differentials determine the associated cone angles. The loci of differentials with prescribed numbers and orders of zeros stratify the moduli space of differentials, and the adjacency relations among strata arise from merging zeros. The action of affine transformations on a translation surface preserves the cone angles of its singularities, thereby inducing a ${\rm GL}_2^{+}(\mathbb R)$-action on each stratum of differentials with a fixed zero type, commonly referred to as Teichm\"uller dynamics. Up to scalar multiplication, holomorphic differentials correspond to canonical divisors, which govern the intrinsic geometry of the underlying algebraic curves. These interrelated yet differently motivated viewpoints make the study of differentials and their strata a central topic at the intersection of several areas of mathematics. 

Over the past few decades, remarkable progress on differentials and their moduli spaces has been achieved in diverse directions, including the study of dynamical invariants on translation surfaces, the classification of ${\rm GL}_2^{+}(\mathbb R)$-orbit closures, compactifications of strata of differentials, intersection number and cycle class computations, and developments in birational geometry. There already exist many excellent surveys and lecture notes on these topics; we simply refer the reader to \cite{Zo06, Mo13, Wr15, Ch17, Ma18, Fi24, AM24} among many others.  

Nevertheless, our understanding of the topological properties of strata of differentials remains quite limited. In this survey, we focus on the topology of strata of differentials from the viewpoint of algebraic geometry. At first glance, it may seem presumptuous to approach topological questions using algebraic-geometric methods. However, many topological invariants of strata of differentials admit natural counterparts in algebraic geometry. Furthermore, standard techniques in algebraic geometry---such as compactifications of moduli spaces, deformations of singularities, and the construction and positivity of characteristic classes---offer powerful tools for exploring topology-related phenomena. 

Concretely, we will discuss the following topics: connected components of strata of differentials, where additional components are distinguished by hyperelliptic and spin structures; the orbifold fundamental group and $K(\pi,1)$-type questions for strata of differentials, which exhibit deep connections with framed mapping class groups, deformations of singularities, and stability conditions; Euler characteristics of strata of differentials, which can be obtained via intersection theory on compactified strata; and tautological characteristic classes on strata of differentials, whose positivity properties yield insights into the affine geometry, boundary structures, and cohomological dimensions of the strata. In addition to holomorphic differentials, we will also address analogous questions for strata of meromorphic differentials, $k$-differentials, and generalized strata with residue constraints.    

Since the topics covered in this survey lie at the intersection of several fields, we aim to prioritize geometric intuition behind the relevant concepts and ideas, while minimizing technical details without sacrificing precision, in order to maintain a natural and coherent exposition. This survey is by no means exhaustive, and omissions or biases are inevitable. Nevertheless, we hope that readers from diverse backgrounds will find it accessible and informative, even on a first or quick reading. 

\subsection*{Acknowledgements} Part of this survey is based on the talk delivered by the first-named author at the session ``Geometry of Moduli'' during the $2025$ Summer Research Institute in Algebraic Geometry. The authors are grateful to the conference and session organizers for the opportunity to present their work and to contribute to the conference proceedings. 

\section{Preliminaries on strata of differentials}
\label{sec:prelim}

In this section, we provide a brief review of the fundamental properties of differentials and their associated moduli spaces; see~\cite[Section 3]{Zo06} for a more comprehensive introduction. 

Let $\mu = (m_1, \ldots, m_n)$ be an integral partition of $2g - 2$, where $m_i \ge 0$ for all $i$. Given a nontrivial holomorphic differential $\omega$ on a smooth Riemann surface $X$ of genus $g$ (also called an Abelian differential), suppose that $\omega$ has $n$ distinct zeros at $z_1, \ldots, z_n$, where each $z_i$ is a zero of order $m_i$. Then we say that $\omega$ is of type $\mu$.

There is a fascinating interplay between holomorphic differentials and translation structures with conical singularities. Away from the zeros of $\omega$, we can locally choose an analytic coordinate $u$ such that $\omega = du$. Viewing $u = x + {\rm i}y$ as the standard coordinate on the Euclidean plane, $\omega$ induces a flat metric on $X \setminus \{z_1, \ldots, z_n\}$. If $v$ is another local coordinate satisfying $\omega = dv$, then $u$ and $v$ differ by a translation, which preserves the induced flat metric. 

At a zero $z_i$ of order $m_i$, one can find a local coordinate $w$ such that $\omega = w^{m_i}dw$. Comparing this with the nearby expression $\omega = du$, we obtain $w = u^{m_i+1}$ up to a scalar factor. Consequently, as one travels once around $z_i$, the total angle under the flat metric induced by $du$ is $2\pi(m_i + 1)$. In this sense, $\omega$ induces a translation structure on $X$ with cone points at $z_1, \ldots, z_n$, where the cone angle at each $z_i$ equals $2\pi(m_i + 1)$.

Let $\mathcal H(\mu)$ be the moduli space parametrizing differentials $(X,\omega)$ of type $\mu$. Locally at $(X, \omega)$, the space $\mathcal H(\mu)$ admits a coordinate system given by the relative cohomology $H^1(X, \{z_1,\ldots, z_n\};\mathbb Z)$. More precisely, suppose 
$$\gamma_1,\ldots, \gamma_{2g}; \gamma_{2g+1},\ldots, \gamma_{2g+n-1}$$ 
form a basis for the relative homology $H_1(X, \{z_1,\ldots, z_n\};\mathbb Z)$, where the first $2g$ paths are closed and span the absolute homology $H_1(X; \mathbb Z)$, and each of the remaining paths $\gamma_{2g+i}$ joins $z_i$ to $z_n$ for $i = 1,\ldots, n-1$. Then 
$$\int_{\gamma_1}\omega,\ldots, \int_{\gamma_{2g}}\omega; \int_{\gamma_{2g+1}}\omega,\ldots, \int_{\gamma_{2g+n-1}}\omega$$ 
provide local coordinates for $\mathcal H(\mu)$, called period coordinates. 

From the viewpoint of the translation surface structure induced by $\omega$, deforming these period coordinates amounts to varying the parallel edges that bound the translation surface while preserving the number and cone angles of the singularities, thus determining a neighborhood of $(X,\omega)$ in $\mathcal H(\mu)$. As a consequence, $\mathcal H(\mu)$ is a complex orbifold of dimension 
$$\dim_{\mathbb C}\mathcal H(m_1,\ldots, m_n) = 2g+n-1,$$ 
where the orbifold structure arises from special $(X,\omega)$ with nontrivial automorphisms. We will also denote by 
$$\mathbb P\mathcal H(\mu) =\mathcal H(\mu) / \mathbb C^{*}  $$
the projectivization, which parametrizes canonical divisors of type $\mu$. 

Note that for fixed $g$, the union of $\mathcal H(\mu)$ over all partitions $\mu$ of $2g-2$ is the Hodge bundle $\mathcal H_g$ (with the zero section removed) of holomorphic differentials over the moduli space $\mathcal M_g$ of genus-$g$ Riemann surfaces. Shrinking a short period joining two zeros of orders $m_i$ and $m_j$ causes them to merge, which induces a natural stratification of $\mathcal H_g$ by the strata $\mathcal H(\mu)$, where $\mathcal H(\mu)$ contains $\mathcal H(\mu')$ in its boundary if and only if $\mu$ can specialize to $\mu'$ by merging its parts. In this sense, $\mathcal H(\mu)$ is called the stratum of holomorphic differentials of type $\mu$. 

From the perspective of dynamics on Riemann surfaces, a key feature of holomorphic differentials is that ${\rm GL}_2^{+}(\mathbb R)$ acts on them by affine transformations, preserving the orders of their zeros. Consequently, this induces a natural ${\rm GL}_2^{+}(\mathbb R)$-action on each stratum of holomorphic differentials with prescribed zero orders, whose study is commonly referred to as Teichm\"uller dynamics. Classfying the orbit closures of 
the ${\rm GL}_2^{+}(\mathbb R)$-action in the period coordinates of the strata is a central problem in this field. A breakthrough result of Eskin, Mirzakhani, and Mohammadi \cite{EMM15, EM18} establishes that every ${\rm GL}_2^{+}(\mathbb R)$-orbit closure is locally defined by linear equations in the period coordinates, whose coefficients were later shown by Filip \cite{Fi16} to be algebraic numbers. Although ${\rm GL}_2^{+}(\mathbb R)$-orbit closures exhibit geometric features similar to those of the ambient strata, their classification and related topological properties remain far less understood. We refer the reader to \cite{Wr16, Fi24} for a more detailed introduction and further references on this topic. 

Besides holomorphic differentials, one may also consider strata of meromorphic differentials with prescribed zero and pole orders, as well as strata of $k$-differentials, which are sections of the $k$th power of the canonical bundle. From the viewpoint of flat metrics and translation structures, a meromorphic differential corresponds to a translation surface of infinite area, with the poles interpreted as points at infinity. Likewise, a $k$-differential corresponds to a $\frac{1}{k}$-translation surface, whose transition functions consist of translations and rotations by angles that are integer multiples of $2\pi/k$. We refer the reader to \cite{Bo15, BCGGMk} for more detailed introductions to the strata of meromorphic differentials and $k$-differentials, respectively. 

\section{Connected components of strata of differentials}
\label{sec:components}

Classifying the connected components of a space---equivalently, determining its $\pi_0$---is the most fundamental step in understanding its topology.

In this section, we first review the hyperelliptic and spin structures that serve to distinguish connected components of the strata of differentials. We then discuss the known results and open problems concerning the classification of connected components of the strata of holomorphic differentials, meromorphic differentials, quadratic differentials, and $k$-differentials, as well as generalized strata of differentials with residue constraints. 

\subsection{The hyperelliptic structure}
\label{subsec:hyp}

Let $X$ be a smooth hyperelliptic curve of genus $g \geq 2$, and let $\tau \in {\rm Aut}(X)$ denote the hyperelliptic involution inducing the double cover $X \to X/\tau \cong \mathbb P^1$. 

Any nontrivial holomorphic differential $\omega$ on $X$ has zeros of the form 
$$\sum_{i=1}^s 2k_i w_i + \sum_{j=1}^t l_j (p_j + q_j),$$ 
where each $w_i$ is a Weierstrass point fixed by $\tau$, and each pair $(p_j, q_j)$ consists of points exchanged by $\tau$. In this situation one has $\tau^{*}\omega = -\omega$, which on the level of the associated translation surface corresponds to a rotation by $180^\circ$. 

The locus of such hyperelliptic differentials has dimension $2g-1+t$, whereas the ambient projectivized stratum $\mathbb P\mathcal H(\mu)$ associated to 
 $$\mu = (2k_1, \ldots, 2k_s, l_1, l_1, \ldots, l_t, l_t)$$ 
 has dimension $2g-2 + s + 2t$. In particular, these two dimensions agree if and only if $\mu = (2g-2)$ or $\mu = (g-1, g-1)$; see Section~\ref{subsec:conn-holo} below for a comparison. 

\subsection{The spin structure}
\label{subsec:spin}

Let $\mu = (2k_1, \ldots, 2k_n)$ be a partition of $2g-2$ whose entries are all even. For any $(X, \omega) \in \mathcal H(\mu)$, the line bundle 
$$L = \mathcal O_X(k_1z_1 + \cdots + k_nz_n)$$
satisfies $L^{\otimes 2} \cong K_X$, and is therefore a theta-characteristic. In this setting, a classical result due to Atiyah \cite{At71} and Mumford \cite{Mu71} states that 
$$ \dim H^0(X, L) \pmod{2}$$ 
is a deformation invariant, which we refer to as the spin parity associated to~$\omega$.  

There is an alternative way to determine the spin parity in terms of the translation structure induced by $\omega$. Given a closed path $\gamma$ in $X$ that does not pass through the zeros of $\omega$, denote by ${\rm Ind}_{\gamma}(\omega)$ the winding number of $\gamma$ with respect to the flat metric defined by $\omega$. If $\gamma$ is deformed to pass across a zero of order $m$, its winding number changes by $m$. Consequently, if $\omega$ has only zeros of even order, then ${\rm Ind}_{\gamma}(\omega) \pmod{2}$ depends only on the homology class represented by $\gamma$. In this situation, let $\{(\alpha_i,\beta_i)\}_{i=1}^g$ be a symplectic basis of $H_1(X;\mathbb Z)$, chosen so that none of the curves intersect the zeros of $\omega$. One then defines the Arf invariant by 
$${\rm Arf}(\omega) = \sum_{i=1}^g ({\rm Ind}_{\alpha_i}(\omega)+1)({\rm Ind}_{\beta_i}(\omega)+1) \pmod{2}.$$
Johnson \cite{Jo80} proved that the Arf invariant agrees with the spin parity, i.e.,  
\begin{align}
\label{eq:spin}
    {\rm Arf}(\omega) = \dim H^0(X, L) \pmod{2}. 
\end{align}

\subsection{Connected components of strata of holomorphic differentials} 
\label{subsec:conn-holo}

Using the hyperelliptic and spin structures described above, Kontsevich and Zorich obtained a complete classification of the connected components of the strata of holomorphic differentials. In concise form (see \cite[Theorems~1 and~2]{KZ03}), a stratum of holomorphic differentials $\mathcal H(\mu)$ may have up to three connected components, with extra components arising from hyperelliptic or spin structures. More precisely, a hyperelliptic component occurs if and only if $\mu = (2g-2)$ or $\mu = (g-1,g-1)$, while odd and even spin components appear exactly when every entry of $\mu$ is even. 

For example, when $g \ge 4$, the minimal stratum $\mathcal H(2g-2)$ has three connected components:
$\mathcal H(2g-2)^{\rm hyp}$, $\mathcal H(2g-2)^{\rm odd}$, and $\mathcal H(2g-2)^{\rm even}$. Here, $\mathcal H(2g-2)^{\rm hyp}$ parametrizes hyperelliptic differentials whose unique zero lies at a Weierstrass point, while $\mathcal H(2g-2)^{\rm odd}$ (resp. $\mathcal H(2g-2)^{\rm even}$) parametrizes non-hyperelliptic differentials $\omega$ for which the associated spin structure
$\dim H^0(X, (g-1)z) \pmod{2}$ is odd (resp. even). We note that in low genus, some hyperelliptic and spin components coincide; for instance,
$\mathcal H(4)^{\rm hyp} = \mathcal H(4)^{\rm even}$ when $g = 3$.

When hyperelliptic or spin structures are present, it is relatively straightforward to see that they lead to additional connected components in the strata of holomorphic differentials. The difficult part, however, lies in proving that these are the only sources of disconnectedness—that is, establishing that no further hidden invariants can distinguish additional connected components. 

The proof by Kontsevich and Zorich relies on a subtle analytic surgery known as bubbling a handle, which inductively constructs differentials in genus 
$g+1$ from those in genus $g$. A second key ingredient comes from the theory of Rauzy classes associated with interval exchange transformations, which ensures that every connected component in higher genus can be reached through repeated applications of this surgery. Together, these ideas provide a complete classification of the connected components of the strata of holomorphic differentials. 

As of today, no alternative proof of the Kontsevich--Zorich classification is known; in particular, there is currently no algebro-geometric proof that would apply over fields other than $\mathbb C$. A promising direction may lie in the recent developments on the multi-scale compactification of the strata of differentials \cite{BCGGM1, BCGGM2}, whose boundary---analogous to that of the Deligne--Mumford compactification of the moduli space of curves---is built from strata of lower genus or with more specialized singularity types. One may hope to use this compactification to inductively determine the connected components of the strata. However, such an inductive strategy requires knowing that every connected component contains specific types of boundary differentials, such as those arising from bubbling a handle, where Kontsevich and Zorich crucially used the non-algebraic input of Rauzy classes. 

We summarize these considerations in the following question. 

\begin{problem}
Instead of using Rauzy classes, find an algebraic argument to show that every connected component of the strata of holomorphic differentials contains boundary differentials obtained by bubbling a handle. 

Consequently, find an algebraic proof of the Kontsevich--Zorich classification such that it can be generalized to fields other than $\mathbb C$. 
\end{problem}

Finally, we remark that one may also study strata of holomorphic differentials on surfaces with Teichmüller markings, i.e., strata defined inside the Hodge bundle over Teichmüller space. The classification of connected components of such marked strata was carried out by Calderon and Salter \cite{Ca20, CS21}. Unlike the situation for unmarked translation surfaces, many more connected components appear over Teichmüller space. In particular, there are infinitely many hyperelliptic connected components for $g \geq 3$; and for $g \geq 5$, whenever a common divisor $r$ divides all zero orders of the differential, the $r$-spin structure (generalizing the case $r=2$) distinguishes exactly $r^{2g}$ non-hyperelliptic components; see also~\cite{Ha21} for some related results. This perspective is closely related to the study of orbifold fundamental groups of strata and their monodromy images in the mapping class group; see Section~\ref{subsec:fundamental}. 

\subsection{Connected components of strata of meromorphic differentials} 

One can also consider the strata of meromorphic differentials with prescribed zero and pole orders. Their connected components were classified by Boissy \cite{Bo15}, using hyperelliptic and spin structures together with inductive procedures similar to those in the work of Kontsevich and Zorich. Nevertheless, two new features arise in the case of meromorphic differentials. 

First, suppose $\mu = (m_1, \ldots, m_n)$ is a signature of meromorphic differentials in genus one, and denote by $\gcd(\mu) = \gcd(m_1, \ldots, m_n)$. Given a meromorphic differential $\omega$ of type $\mu$ on an elliptic curve $E$, define the rotation number 
$$ {\rm rot} (\omega) = \gcd ({\rm Ind}_{\alpha}(\omega), {\rm Ind}_{\beta}(\omega), m_1, \ldots, m_n), $$
where $\alpha$ and $\beta$ form a symplectic basis of $H_1(E; \mathbb Z)$, and ${\rm Ind}$ denotes the winding number of a closed path in $E$ with respect to the flat metric induced by $\omega$. Boissy showed that these rotation numbers distinguish the connected components of the strata of genus-one meromorphic differentials; see~\cite[Theorem 4.3]{Bo15}. 

One can also interpret rotation numbers from the viewpoint of algebraic geometry. Let $m_1z_1 + \cdots + m_n z_n \sim K_E \cong \mathcal O_E$ be the meromorphic canonical divisor associated to $\omega$. Let $d$ be the largest divisor of $\gcd(\mu)$ such that  
$$(m_1/d)z_1 + \cdots + (m_n/d)z_n \sim \mathcal O_E,$$ 
where $d$ is called the torsion number of $\omega$. It is not hard to see that torsion numbers can also be used to distinguish connected components of the stratum of genus-one meromorphic differentials of type $\mu$; see, e.g.,~\cite[Section 3.2]{CC14}. In fact, the rotation number and torsion number of $\omega$ coincide; see~\cite[Section 3.4]{Ta18} and \cite[Proposition 3.13]{CG22}. 

Next, poles of a meromorphic differential are located at points at infinity with respect to the corresponding translation surface structure. Under the vertical geodesic flow, meromorphic differentials (without vertical saddle connections) admit a decomposition consisting of (broken) half-infinite cylinders (in the case of simple poles) and (broken) half-planes (in the case of higher-order poles), where the boundaries of these cylinders and planes consist of concatenations of saddle connections with positive real parts. This allows Boissy to use a more geometric argument, instead of Rauzy classes, by shrinking certain short boundary saddle connections to reach every connected component in higher genus through bubbling a handle from lower genus; see~\cite[Proposition 6.1]{Bo15}. 

\subsection{Connected components of strata of quadratic differentials} 

A quadratic differential $q$ on a Riemann surface $X$ is a section of the bi-canonical bundle $K_X^{\otimes 2}$. Taking a local square root of $q$ induces a half-translation surface structure on $X$, where the transition functions are translations and rotations by angles that are multiples of $\pi$. In addition, there exists a canonical double cover $\pi\colon \widetilde{X}\to X$ such that
$\pi^{*}q = \widetilde{\omega}^2$, where $\pi$ is branched at the odd singularities of $q$, $\widetilde{\omega}$ is a globally defined Abelian differential on $\widetilde{X}$, and $\widetilde{\omega}$ is anti-invariant under the involution of the canonical double cover. Note that $\widetilde{\omega}$ is holomorphic if and only if $q$ has at worst simple poles, where a simple pole of $q$ corresponds to a conical singularity of angle $\pi$ in the associated half-translation surface. 

Many properties of Abelian differentials extend naturally to quadratic differentials. In particular, the ${\rm GL}_2^{+}(\mathbb R)$-action preserves the orders of singularities of quadratic differentials and therefore acts on each stratum of quadratic differentials with prescribed zero and pole orders. Moreover, the image of a stratum of quadratic differentials lifted via the canonical double cover into the associated stratum of Abelian differentials is an invariant subspace under the ${\rm GL}_2^{+}(\mathbb R)$-action. 

Denote by $\mathcal Q(\mu)$ the stratum of primitive quadratic differentials of type $\mu$, i.e., those that are not global squares of Abelian differentials. When all entries of $\mu$ are greater than or equal to $-1$---equivalently, when the corresponding half-translation surfaces have finite area---the connected components of $\mathcal Q(\mu)$ were classified by Lanneau \cite{La08}, using inductive procedures similar to those in \cite{KZ03}. When $\mu$ contains higher-order poles, the classification of connected components of $\mathcal Q(\mu)$ was completed in \cite{CG22}, using ideas parallel to those in \cite{Bo15}. As in the case of meromorphic differentials, several new features arise for quadratic differentials. 

First, the hyperelliptic structure, when it exists, can still be used to distinguish connected components of $\mathcal Q(\mu)$; see \cite{La04H}. However, the spin structure does not produce additional components in this setting; see \cite{La04S}. More precisely, if every entry of $\mu$ is not $2 \pmod{4}$, then $\mathcal Q(\mu)$ lifts, via the canonical double cover, to a stratum of Abelian differentials whose zero and pole orders are all even, and hence admits both odd and even spin components. Remarkably, it was shown in \cite[Theorems 1.1 and 1.2]{La04S} that exactly one of these two spin components contains the lifted image of $\mathcal Q(\mu)$. In other words, the spin parity on the canonical double cover is completely determined by the signature $\mu$ of the quadratic differential. To prove this result, Lanneau constructed an explicit symplectic basis for the canonical double cover and computed the associated Arf invariant. In light of the relation~\eqref{eq:spin}, it would be interesting to recover this result by analyzing the dimension of the space of global sections of the corresponding theta-characteristic line bundle. 

Next, in genus three and four, there are several exceptional strata of primitive quadratic differentials:  
\begin{align*}
g = 3\colon & \quad \mathcal Q(9, -1), \quad \mathcal Q(6,3, -1), \quad \mathcal Q(3,3,3, -1); \\
g = 4\colon & \quad \mathcal Q(12), \quad \mathcal Q(9,3), \quad \mathcal Q(6, 6), \quad \mathcal Q(6, 3,3), \quad \mathcal Q(3,3,3,3).
\end{align*}
Each of these strata consists of two connected components, referred to as the ``regular'' and ``irregular'' components, which cannot be distinguished by the hyperelliptic or spin structures. Instead, they are distinguished using the Brill--Noether geometry associated to one third of the zero divisors of the corresponding quadratic differentials; see~\cite[Sections 6 and 7]{CM14}. 

For instance, let $9z - p$ be the bi-canonical divisor associated to a quadratic differential $q$ on a genus-three curve $X$ parametrized in $\mathcal Q(9,-1)$. Then $(X, q)\in \mathcal Q(9,-1)^{\rm reg}$ if $h^0(X, 3z) = 1$, and $(X, q)\in \mathcal Q(9,-1)^{\rm irr}$ if $h^0(X, 3z) = 2$. The proof of this result makes use of the canonical embedding of $X$ into $\mathbb P^2$ as a plane quartic, which intersects a special plane cubic $E$ whose associated $3$-torsion geometry distinguishes the two cases. However, no flat-geometric invariant is currently known that distinguishes the regular and irregular components. 

We summarize these discussions in the following question. 

\begin{problem}
For a stratum of quadratic differentials of spin type---i.e., when the square-root Abelian differential on the canonical double cover has only even-order singularities---determine the spin parity via the dimension of the space of global sections of the associated theta-characteristic line bundle.

Additionally, for the regular and irregular connected components of the exceptional strata of quadratic differentials in genera three and four, find a flat-geometric invariant capable of distinguishing between them. 
\end{problem}

\subsection{Connected components of strata of $k$-differentials} 

In general, one can also study the strata $\mathcal H^k(\mu)$ of $k$-differentials with zero and pole orders specified by $\mu$; see~\cite{BCGGMk} for an introduction to $k$-differentials. However, for $k \geq 3$, very little is known about the connected components of $\mathcal H^k(\mu)$. Among the few existing results, see~\cite[Theorems 1.1 and 1.2]{CG22} for a discussion of hyperelliptic and spin structures for $k$-differentials; see~\cite[Corollary 1.6]{BG25} for a classification of connected components of the minimal strata $\mathcal H^k(2k)$ in genus two; and see Section~\ref{subsec:simple} for the case when there are sufficiently many simple zeros. 

One obstruction to studying $k$-differentials for $k \geq 3$ is the lack of a ${\rm GL}_2^{+}(\mathbb R)$-action. As a consequence, even in the presence of higher-order poles---that is, when the corresponding $\tfrac{1}{k}$-translation surface has infinite area---we do not have a clean decomposition of the polar domains, due to the absence of a well-defined vertical geodesic flow. This contrasts with the case of meromorphic differentials treated in \cite{Bo15}, where such a decomposition plays a key role. In particular, shrinking a short saddle connection on the boundary of a polar domain may cause the surface to degenerate; see~\cite[Remark 6.8]{CG22} for a detailed discussion of this issue. 

We summarize these considerations in the following question. 

\begin{problem}
    Find new invariants and structures, if they exist, for $k$-differentials, and classify the connected components of the strata of $k$-differentials for all $k$.  
\end{problem}

\subsection{Connected components of generalized strata of differentials} 
\label{subsec:generalized}

In the multi-scale compactification of strata of differentials, a key mechanism used to isolate the main component of smoothable degenerate differentials is the global residue condition; see~\cite[Section 2.4 (4)]{BCGGM2}. Roughly speaking, when differentials on Riemann surfaces degenerate to a nodal differential, if the vanishing cycles of certain nodes sum to the trivial homology class, then by Stokes' theorem one obtains a sum-to-zero constraint on the residues of the corresponding meromorphic differentials at those nodes. These meromorphic differentials occur on subsurface components that collapse in the degeneration but reappear after expanding the complementary components to infinity; see~\cite[Figures 15 and 16]{Ch17JDG} for illustrative examples. 

Consequently, it is natural to study loci of meromorphic differentials satisfying prescribed residue relations---these are referred to as generalized strata of differentials. In fact, even when one is ultimately interested only in strata of differentials without residue constraints, such generalized strata inevitably appear when inductively analyzing boundary structures via the multi-scale compactification, due to the global residue condition, and therefore must be taken into account. 

The most special case of generalized strata of differentials is when the residue at every pole is required to be zero; these are called strata of residueless differentials. The classification of connected components of strata of residueless differentials was carried out by Lee~\cite{Lee24}, where hyperelliptic structure, spin parity, and rotation number (in genus one) again play the major roles. More generally, for generalized strata of differentials in which certain specified collections of residues sum to zero, the classification of connected components was completed by Lee and Wong~\cite{LW25}. 

Counterintuitively, the only remaining case not covered in these works is the strata of residueless differentials in genus zero. This case is closely related to the longstanding Hurwitz realization problem for branched covers of the sphere with prescribed ramification profiles, where the integral of a residueless differential yields the rational function associated to the corresponding branched cover. 

One can also define residues for poles of $k$-differentials; see~\cite[Section 3.1]{BCGGMk}. A related fundamental question concerns the existence of $k$-differentials with prescribed singularity orders and specified residues at each pole. This problem was studied in a series of works by Gendron and Tahar~\cite{GT21, GT25, GTk25}. Once again, the only major missing case is the realization problem for residueless differentials in genus zero with prescribed zero and pole orders. 

We summarize these discussions in the following question. 

\begin{problem}
    Determine when strata of residueless differentials in genus zero with prescribed zero and pole orders are nonempty, and classify their connected components. 
    
    Consequently, resolve the Hurwitz realization problem for branched covers of $\mathbb P^1$ with prescribed ramification profiles.  
\end{problem}

\section{Fundamental groups of strata of differentials}
\label{sec:fundamental}

For the remainder of the survey, as a convention, when referring to a stratum or using notation such as $\mathcal H(\mu)$, we will mean a single connected component in cases where the stratum is disconnected.

As mentioned previously, beyond the classification of connected components, most topological properties of the strata of differentials remain poorly understood. In this section, we review recent advances toward understanding fundamental groups of strata of differentials, as well as connections to other topics such as 
the monodromy representation, deformations of singularities, and stability conditions. 

\subsection{The $K(\pi,1)$ conjecture}
\label{subsec:fundamental}

It is a classical fact that the moduli space of pointed Riemann surfaces ${\mathcal M}_{g,n}$ is the quotient of the Teichm\"uller space ${\mathcal T}_{g,n}$ by the mapping class group ${\rm Mod}_g^n$, where ${\mathcal T}_{g,n}$, as the orbifold universal cover of ${\mathcal M}_{g,n}$, is contractible, and the orbifold fundamental group $\pi_1^{\rm orb}({\mathcal M}_{g,n})$ can be identified with the mapping class group ${\rm Mod}_g^n$.

Kontsevich and Zorich~\cite[Section 4]{KZ97} conjectured that the analogous statement holds for the strata of holomorphic differentials, known as the $K(\pi,1)$ conjecture for strata: namely, that the orbifold universal cover of $\mathcal H(\mu)$ is contractible and that the orbifold fundamental group $\pi_1^{\rm orb}(\mathcal H(\mu))$ is commensurable with a suitable mapping class group.

This $K(\pi,1)$ conjecture can be verified for the hyperelliptic strata, essentially because the locus of hyperelliptic curves equipped with canonical divisors of fixed type can be identified---up to a finite base change---with the configuration space ${\mathcal M}_{0,n}$ of $n$ marked points on the sphere. Here, the $n$ points correspond to the branch points of the hyperelliptic double cover along with the images of conjugate pairs appearing in the corresponding canonical divisors, and the base change amounts to labeling these points when necessary. In particular, the mapping class groups in the $K(\pi,1)$ conjecture for hyperelliptic strata are certain braid groups. We refer to \cite[Section 1.4]{LM14} for a detailed discussion of this case.

In the same paper, Looijenga and Mondello also obtained presentations for the orbifold fundamental groups of several non-hyperelliptic strata in genus three, using the geometry of the associated canonical curves and root-system data; however, their results do not yet lead to a conclusive verification of the $K(\pi,1)$ conjecture for these strata. 

\subsection{The monodromy representation}
\label{subsec:monodromy}

A closely related study is about the topological monodromy representation 
\begin{align}
\label{eq:monodromy}
\rho\colon \pi_1^{\rm orb}(\mathcal H(\mu))\to {\rm Mod}_g^n, 
\end{align}
where a loop in $\mathcal H(\mu)$ based at a differential $(X, \omega)$ 
induces a self-homeomorphism of $X$ preserving the zeros of $\omega$. Calderon and Salter \cite{CS23} studied the monodromy representation $\rho$ and showed that for non-hyperelliptic strata of $g\geq 5$, it surjects onto the subgroup ${\rm Mod}_g^n[\bar{\phi}]$ of the mapping class group which preserves a fixed framing $\bar{\phi}$ of $X$, where $\bar{\phi}$ is induced by the horizontal vector field of $1/\omega$. 

However, knowing the image of the monodromy representation $\rho$ is not sufficient to fully determine the orbifold fundamental group of the strata, since $\rho$ may have a large kernel even in low genus. Indeed, Giannini \cite{Gi24, Gi25} showed that for $\mu = (4)^{\rm nonhyp}$, $(3,1)$, and $(6)^{\rm even}$, the kernel of $\rho$ contains a non-abelian free group of rank $2$. The proof of these results builds on Wajnryb’s work \cite{Wa99} and involves an analysis of Artin groups associated with the Dynkin diagrams of singularities of $ADE$ type. 

We remark that one can also study the fundamental groups of strata from a more dynamical perspective. For example, Bell, Delecroix, Gadre, Guti\'errez-Romo, and Saul Schleimer \cite{BDGGS} showed that the fundamental group of a stratum component of rooted quadratic differentials is generated by loops formed by the diagonal flow. 

\subsection{Deformations of singularities}
\label{subsec:singularities}

To see why such singularities arise in the study of differentials, consider the curve singularity of type $A_{2g}$ defined by 
$$y^{2} - x^{2g+1} = 0.$$ 
Its miniversal deformation is parametrized by  
$$ y^2 - x^{2g+1} + c_{2g-1} x^{2g-1} + \cdots + c_{1} x + c_0 = 0. $$
In particular, every smooth deformation is a hyperelliptic curve $X$ of genus $g$, where the hyperelliptic double cover is given by $(x, y) \mapsto x$. Moreover, $X$ carries a holomorphic differential $\omega$ with a unique zero $z$ at infinity on the affine chart. Therefore, we may identify the hyperelliptic stratum $\mathcal H(2g-2)^{\rm hyp}$ with the locus of smooth deformations of the $A_{2g}$ singularity. 

In this sense, beyond the $A_{2g}$ singularity, the classical singularities of types $A_{2g+1}$, $D_{2g+1}$, $D_{2g+2}$, $E_6$, $E_7$, and $E_8$ correspond, respectively, to the strata $\mathcal H(g-1,g-1)^{\rm hyp}$, $\mathcal H(2g-2, 0)^{\rm hyp}$, $\mathcal H(g-1,g-1, 0)^{\rm hyp}$, $\mathcal H(4)^{\rm odd}$, 
$\mathcal H(3,1)$, and $\mathcal H(6)^{\rm even}$, where an entry $0$ corresponds to an ordinary marked point, i.e., a ``zero'' of order $0$. It is worth noting that all these singularities are Gorenstein with a natural $\mathbb G_m$-action, where Gorenstein means that the dualizing sheaf of the underlying singular curve remains a holomorphic line bundle, and the $\mathbb G_m$-action arises from assigning weights to the parameters in the defining equations of the singularities---for instance, giving $x$ and $y$ weights $2$ and $2g+1$, respectively, for the $A_{2g}$ singularity defined by $y^2 - x^{2g+1} = 0$. 

In general, Chen and Yu \cite{CY25} used deformations of isolated Gorenstein curve singularities with a $\mathbb G_m$-action to study strata of holomorphic differentials. Beyond the $ADE$ singularities, several new singularities appear in low genus, each corresponding to an entire stratum (or connected component) in which every differential $(X, \omega)$ admits an isotrivial degeneration (with respect to the $\mathbb G_m$-action) to the same curve singularity. Such strata are called nonvarying, originally discovered and studied due to the uniform behavior of Teichmüller curves contained in them; see~\cite{CM12, YZ13, CM14}. We also refer to Section~\ref{sec:tauto} for the special property that these nonvarying strata are affine varieties with trivial tautological rings. However, for $g \geq 6$, this nonvarying phenomenon disappears (except in the hyperelliptic strata), indicating a subtle hierarchical structure among those singularities whose deformation spaces correspond to the same stratum of differentials.

Returning to the original $K(\pi,1)$ question for strata of differentials, there is a parallel story in the study of deformations of singularities. Roughly speaking, for certain classes of singularities (e.g., positive-dimensional isolated complete intersection singularities), it is conjectured that the discriminant complement---namely, the locus of smooth deformations---in the miniversal deformation space satisfies the $K(\pi,1)$ property. From the perspective of associated Artin groups, this conjecture is sometimes attributed to Arnold, Brieskorn, Pham, and Thom; see \cite[p.~26]{Ar88}, \cite[p.~185]{Lo84}, and \cite{Pa14} for related discussions. 
The conjecture was verified for $ADE$ singularities by Brieskorn and Deligne \cite{Br73, De72}, mirroring the situation for the corresponding hyperelliptic strata and several nonvarying strata in low genus. However, the conjecture remains open for many other singularities, including those corresponding to the remaining nonvarying strata classified in \cite[Theorem 1.2]{CY25}. 

We summarize these considerations in the following question. 

\begin{problem}
    For the nonvarying strata of holomorphic differentials that are not of $ADE$ type, study the corresponding monodromy kernels and deformation spaces of the associated singularities, and determine whether the $K(\pi,1)$ property holds in these cases. 
   
    For varying strata, classify the associated Gorenstein singularities, describe how their deformation spaces assemble to form the stratum, and use this structure to investigate the topology of the full stratum.
\end{problem}

\subsection{Differentials with many simple zeros}
\label{subsec:simple}

There is one regime of signatures $\mu$ where the fundamental group and the monodromy representation of $\mathcal H(\mu)$ can be understood---namely, when there are sufficiently many simple zeros. 

Let $\mu = (1^{n-k}, m_1, \ldots, m_k)$ where $m_i > 1$ for the last $k$ zeros. Associating $(X, \omega) \in \mathcal H(\mu)$ to the underlying curve $X$ marked at the higher-order zeros $z_1, \ldots, z_k$ induces a map 
$$\mathcal H(\mu)\to \mathcal M_{g,k},$$ 
which exhibits $\mathcal H(\mu)$ as an open subset of an almost “vector bundle--like’’ structure over $\mathcal M_{g,k}$. The completion of the fiber over $(X, z_1, \ldots, z_k)$ is given by 
$$H^0(X, K - m_1z_1 - \cdots - m_kz_k)$$ and the complement of $\mathcal H(\mu)$ in each fiber is the discriminant hypersurface where two or more zeros coincide. 

This is not quite a vector bundle because, for special $(X, z_1,\ldots, z_k)$---for instance, those with Brill–Noether special geometry---the dimension of $H^0(X, K - m_1z_1 - \cdots - m_kz_k)$ may jump. However, the degree of the divisor class $K - m_1z_1 - \cdots - m_kz_k$ is $n-k$, which equals the number of simple zeros in $\mu$. Therefore, when there are sufficiently many simple zeros, the locus of such special $(X, \omega)$ has high codimension in $\mathcal M_{g,k}$. In particular, removing these loci does not affect the fundamental group nor the (co)homology groups in low (co)dimension; see~\cite[Examples 2.3--2.5]{Ch19} for some related examples and discussions. 

Salter adapted this idea to show that the monodromy representation $\rho$ in~\eqref{eq:monodromy} is injective when there are sufficiently many simple zeros. In this range, the orbifold fundamental group of the stratum can therefore be identified with the framed mapping class group; see~\cite[Theorem A]{Sa25} for a precise bound on the number of required simple zeros. 

Salter’s argument builds on work of Shimada~\cite{Sh10}, who studied fundamental groups of algebraic varieties that behave like fiber bundles away from loci of high codimension. A key example is the fundamental group of the locus of reduced divisors $|L|_{\rm red}$ in the linear system $|L|$ for a very ample line bundle $L$ on a smooth curve $X$, where $L$ varies in the Picard group of $X$ and induces a fiber bundle--like structure via the associated linear systems. 

An additional difficulty arises from the fact that, for signatures such as $\mu = (1^{2g-2})$, one can study the fundamental group of $|L|_{\rm red}$ for a general line bundle $L$ of degree $2g-2$ in the Picard group using the Abel–Jacobi map 
$$X^{2g-2} \to {\rm Pic}^{2g-2}(X).$$ 
However, the canonical bundle $K$ is a highly special point of ${\rm Pic}^{2g-2}(X)$, and its sections correspond precisely to the differentials parametrized in the strata $\mathcal H(\mu)$. To bridge this gap, Salter employed what is referred to as the Shimada--Severi method. The key idea is to use the connectedness of the Severi variety of plane curves to construct an equisingular deformation between generic plane sections of $|L|_{\rm red}$ and $|K|_{\rm red}$, assuming that $X$ has sufficiently generic Weierstrass points. This deformation then induces an isomorphism between their fundamental groups; see~\cite[Section 5]{Sa25}. 

Finally, we remark that although this statement does not explicitly appear in \cite{Sa25}, the same method can be adapted to treat the strata of $k$-differentials with sufficiently many simple zeros. 

\subsection{Differentials and stability conditions}
\label{subsec:stability}

An analogous result has also been established for the fundamental groups of the principal strata of quadratic differentials. When the quadratic differentials carry sufficiently many simple zeros, the associated Abel--Jacobi map and the image of the monodromy representation were first analyzed by Walker \cite{Wa10} using classical methods. More recently, a complete understanding of the fundamental group---namely, a description of the monodromy kernel---was achieved by King and Qiu \cite{KQ20} for the strata of quadratic differentials with simple zeros and higher-order poles, using Bridgeland stability conditions and categorical techniques. This work was subsequently extended by Qiu \cite{Qi25} to include strata of quadratic differentials with simple zeros, higher-order poles, and also double poles. 

The connection between differentials and stability conditions is an exciting and rapidly developing research direction: several striking results have already been obtained, yet many fundamental questions remain open. Roughly speaking, moduli spaces of framed differentials correspond to stability manifolds, where a differential plays the role of a stability condition by assigning periods to framed homology classes, thereby determining the central charge in the sense of Bridgeland stability. These stability manifolds admit a wall-and-chamber decomposition in which the (real codimension-one) walls correspond to loci where two non-homologous classes acquire collinear periods. Moreover, (semi)stable objects correspond to certain geodesics---such as saddle connections joining zeros or closed geodesics filling a cylinder---on the flat surface defined by the differential. The enumeration of these objects is closely tied to BPS states and Donaldson--Thomas invariants and governed by wall-crossing phenomena.

This connection between differentials and stability conditions was first conjectured by Kontsevich and Soibelman \cite[Section 1.5 4)]{KS08} and further motivated by the GMN differentials studied by Gaiotto, Moore, and Neitzke \cite{GMN13}. In the case of meromorphic quadratic differentials with simple zeros, the correspondence was established in full detail by Bridgeland and Smith \cite{BS15}. Subsequently, Haiden, Katzarkov, and Kontsevich \cite{HKK17} extended the framework using Fukaya categories of punctured surfaces, providing a conceptual description of the correspondence for all strata of quadratic differentials. 

Returning to the $K(\pi,1)$ question, we note that---analogous to the situation for strata of differentials and for discriminant complements in deformation spaces of singularities---there is a parallel speculation concerning the contractibility of stability manifolds. Under the correspondence between differentials and stability conditions, stability manifolds (understood as strata of framed differentials) correspond to the universal covers of the strata of unmarked differentials. In this sense, the contractibility problem for stability manifolds and the $K(\pi,1)$ conjecture for strata of differentials become equivalent formulations of the same question. 

We summarize these considerations in the following question. 

\begin{problem}
    Using the correspondence between differentials and stability conditions, study the topology of strata of differentials and the associated stability manifolds from both perspectives. 
\end{problem}

\section{Euler characteristics of strata of differentials}
\label{sec:euler}

As discussed in Section~\ref{sec:fundamental}, even though the fundamental groups and low-degree (co)homology groups of strata of differentials can be determined in certain cases---such as low genus or when many simple zeros are present---the general topological properties of these strata remain largely unknown. Nonetheless, even without full knowledge of the Betti numbers, one can still compute Euler characteristics through alternative methods. In particular, recent progress has been made using intersection theory on compactified strata, and in some cases for loci of meromorphic differentials with prescribed residues. In this section, we review these developments. 

\subsection{The cotangent bundle of the multi-scale compactification}
\label{subsec:cotangent}

Given a smooth complex variety $M$, it is a classical fact that the Euler characteristic $\chi(M)$ can be computed, up to sign, as the degree of the top Chern class of the cotangent bundle $\Omega_M^1$. When $M$ is an orbifold equipped with a normal crossing boundary divisor $D$, the same principle applies: the orbifold Euler characteristic of the interior, $\chi(M \setminus D)$, can be determined by the degree of the top Chern class of the logarithmic cotangent bundle $\Omega_M^1(\log D)$. 

Costantini, M\"oller, and Zachhuber applied this idea to the multi-scale compactification of strata of Abelian differentials and obtained an intersection-theoretic formula for the orbifold Euler characteristics of the strata; see~\cite[Theorem 1.3]{CMZ22}. Implementing this formula requires knowledge of the top self-intersection number of the tautological line bundle class $c_1(\mathcal O(-1))$ on every boundary stratum, which in particular involves generalized strata with residue constraints imposed by the global residue condition, as discussed in Section~\ref{subsec:generalized}. In principle, these intersection numbers can be computed recursively using divisor class relations on the strata. However, the resulting calculations can be highly intricate and, even in low genus, may require substantial computational assistance; see~\cite{CMM23}.

Several follow-up works are now available. The same ideas and techniques can be applied to compute the Euler characteristics of linear subvarieties in the strata which are locally defined by linear equations of the period coordinates, 
such as strata of $k$-differentials lifted into strata of Abelian differentials via the canonical cyclic $k$-cover; see~\cite{CMS25}. Further refinements for the minimal strata and for the spin components appear in \cite{CSS25, Wo24}. 

Finally, we remark that the Euler characteristics of the moduli spaces of curves $\mathcal M_{g,n}$ can also be computed via intersection-theoretic methods; see~\cite{GLN23}. The original proof, due to the seminal work of Harer and Zagier \cite{HZ86}, instead relies on a cell decomposition of $\mathcal M_{g,n}$ based on ribbon graphs, together with graph-counting and related combinatorial techniques. It would be highly desirable to develop analogous cell decompositions for strata of differentials and to derive similar applications. Recent progress in this direction involves studying Delaunay triangulations on strata \cite{Zy22, FZ25}, although a systematic framework and definitive results remain open. 

We summarize these discussions in the following question. 

\begin{problem}
    Using intersection theory or cell decompositions to obtain simplified expressions for the Euler characteristics of strata of differentials, and to uncover additional topological properties of the strata. 
\end{problem}

\subsection{The isoresidual fiberation}
\label{subsec:isoresidual}

As discussed in Section~\ref{subsec:generalized}, meromorphic differentials with residue constraints naturally arise in the boundary of compactified strata. It is therefore meaningful to study loci of meromorphic differentials with prescribed residues at the poles. Let $\mathcal H(a_1,\ldots, a_n; -b_1, \ldots, -b_p)$ be a stratum of meromorphic differentials, where $a_i \ge 0$ and $b_i > 0$. Associating to each differential its residue tuple $\underline{r} = (r_1,\ldots, r_p)$ at the poles induces an isoresidual map
$$ {\rm res}\colon \mathcal H(a_1,\ldots, a_n; -b_1, \ldots, -b_p)\to \mathcal R_p\subset \mathbb C^{p}, $$
whose image lies in the hyperplane $\mathcal R_p$ defined by the relation $r_1 + \cdots + r_p = 0$, in accordance with the Residue Theorem. Consequently, the locus of meromorphic differentials with prescribed residue tuple $\underline{r}$ is precisely the isoresidual fiber $\mathcal F_{\underline{r}}$ of the map~${\rm res}$. 

Note that $$\dim \mathcal H(a_1,\ldots, a_n; -b_1, \ldots, -b_p) = 2g-2 + n + p \geq p-1 = \dim \mathcal R_p,$$ 
where equality holds if and only if $g=0$ and $n=1$, namely for strata of meromorphic differentials on $\mathbb P^1$ with a unique zero. In this case, each isoresidual fiber $\mathcal F_{\underline{r}}$ is zero-dimensional, so studying its topology reduces to determining its cardinality. Using flat-geometric descriptions, Gendron and Tahar computed the cardinality of $\mathcal F_{\underline{r}}$ for generic $\underline{r}$ without partial sum vanishings, i.e., when the sum of the given residue values is nonzero for any proper subset of the poles. 

A complete formula for the cardinality of $\mathcal F_{\underline{r}}$ for $\underline{r}$ with arbitrary partial sum vanishings was obtained by Chen and Prado \cite{CP25}, using intersection-theoretic techniques. The upshot is that associating to a differential the ratio of two residues yields a section of the universal line bundle $\mathcal O(1)$ on the projectivized strata, whose top self-intersection number counts the number of differentials with prescribed residues, with corrections from the boundary of degenerate differentials when partial sum vanishings occur in the residue tuple. Some special cases of this problem, as well as related perspectives via polynomial maps and fixed-point multipliers, can also be found in \cite{BR24, Su17, Su23}. 

Next, one can study when the isoresidual fibers are one-dimensional. There are two cases: differentials on $\mathbb P^1$ with two zeros, or residueless differentials in genus one with a unique zero. For the former, the Euler characteristics of generic one-dimensional isoresidual fibers were determined in \cite{CGPT25}, and for the latter, both the Euler characteristics and the connected components were determined by Lee and Tahar \cite{LT25}. Another variant is to study isoresidual fibers for $k$-differentials and $k$-residues; see~\cite{CGPT-k}. We also remark that the cycle classes of the most special isoresidual fibers---the loci of residueless differentials---satisfy a partial cohomological field theory in $\overline{\mathcal M}_{g,n}$; see~\cite{BPZ24}. The base case of residueless differentials on $\mathbb P^1$ with two zeros appears in the foundational work of Eskin, Masur, and Zorich \cite{EMZ03}, which describes the principal boundary of strata with respect to homologous saddle connections joining the two zeros; see also \cite{CC19} for an interpretation via the Deligne--Mumford compactification. 

Despite these developments, very little is known about the geometry and topology of higher-dimensional isoresidual fibers. We summarize these discussions in the following question. 

\begin{problem}
    Determine the connected components and Euler characteristics of arbitrary isoresidual fibers. Use these results to study the geometry and topology of the strata containing such isoresidual fibers. 
\end{problem}

\section{Tautological rings and affine geometry of strata of differentials}
\label{sec:tauto}

Some research approaches originating from algebraic geometry, although not purely topological, share strong parallels with topology and can be used to derive key topological properties. Examples include comparisons between the Chow ring and the cohomology ring---especially the analogy between tautological classes and characteristic classes arising from vector bundles---as well as the correspondence between affine varieties and Stein manifolds, whose homology groups and boundaries exhibit special structures. In this section, we review recent developments related to the moduli spaces of differentials from these perspectives. 

\subsection{Tautological classes on strata of differentials}
\label{subsec:tauto}

On a moduli space, the characteristic classes of naturally defined vector bundles are commonly referred to as tautological classes. 

For example, consider the moduli space $\mathcal M_{g,n}$ of genus-$g$ curves with $n$ marked points. Let $\lambda_i = c_i(\mathcal H)$ denote the $i$th Chern class of the Hodge bundle $\mathcal H$ of holomorphic differentials on $\mathcal M_{g,n}$. Let $\kappa_i = \pi_{*}\bigl[c_1(\omega_{\pi})^{i+1}\bigr]$ be the Miller--Morita--Mumford classes, where $\pi\colon \mathcal X \to \mathcal M_{g,n}$ is the universal curve and $\omega_{\pi}$ is the relative dualizing bundle. For each marked point, let $\psi_i$ be the first Chern class of the line bundle whose fiber over a pointed curve is the cotangent line at the $i$th marked point. These $\lambda$-, $\kappa$-, and $\psi$-classes generate the tautological ring of $\mathcal M_{g,n}$ (as a subring of either the Chow ring or the cohomology ring). The study of the tautological ring of $\mathcal M_{g,n}$, as well as relations among tautological classes, has been an active and influential topic over the past several decades. To avoid straying too far from the main focus of this survey, we simply refer the interested reader to \cite{Va03} as a starting point.

Given a partition $\mu = (m_1, \ldots, m_n)$ of $2g-2$, denote by $\mathbb P\mathcal H(\mu)$ the stratum of holomorphic differentials of type $\mu$, up to scaling. Equivalently, $\mathbb P\mathcal H(\mu)$ parametrizes canonical divisors of type $\mu$ on smooth genus-$g$ complex algebraic curves. Note that the $\lambda$-, $\kappa$-, and $\psi$-classes pull back to $\mathbb P\mathcal H(\mu)$ from $\mathcal M_{g,n}$. In addition, let 
$$\eta = c_1(\mathcal O(-1)),$$ 
where $\mathcal O(-1)$ is the tautological line bundle associated with the projectivization $\mathbb P\mathcal H(\mu) = \mathcal H(\mu)/\mathbb C^{*}$; that is, the fiber of $\mathcal O(-1)$ over $[(X,\omega)]$ is spanned by the differential $\omega$. Thus, it is natural to define the tautological ring of $\mathbb P\mathcal H(\mu)$ as the subring of its Chow (or cohomology) ring generated by the $\lambda$-, $\kappa$-, $\psi$-, and $\eta$-classes.

On the one hand, compared to $\mathcal M_{g,n}$,
the tautological ring of a stratum of holomorphic differentials has a much simpler structure: it is generated solely by the tautological line bundle class $\eta$, as shown in~\cite{Ch19}. For instance, the relation
\begin{align}
\label{eq:eta-psi}
\eta = (m_i+1)\psi_i
\end{align}
holds for all $i$. This identity can be understood intuitively as follows. At a singularity of order $m_i$, the differential $\omega$ is locally given (up to scaling) by $d(z_i^{m_i+1})$, while $dz_i$ spans the cotangent line at that point, where $z_i$ is a local coordinate. See~\cite[Proposition 2.1]{Ch19} for further explanation, including additional relations involving the remaining tautological classes. These tautological relations also play an important role in the study of the affine geometry of strata of differentials; see Section~\ref{subsec:affine} below. 

On the other hand, it remains largely unknown whether $\eta$ is nonzero in $\mathbb P\mathcal H(\mu)$ and, if so, which powers of $\eta$ remain nonvanishing. For the hyperelliptic strata and the nonvarying strata in low genus, it is known that $\eta$ is trivial, and hence the tautological ring in positive degree is trivial; see~\cite[Theorem 1.1]{Ch24}. In contrast, when there are sufficiently many simple zeros, $\eta$ is nontrivial; see~\cite[Example 2.4]{Ch19}. Indeed, the first-named author conjectured that the rational Picard group (or the second cohomology group) is of rank one generated by $\eta$ for all varying strata of holomorphic differentials; see~\cite[Remarks 3.1 and 3.2]{Ch24} for further discussion. 

As a comparison, the rational Picard group of $\mathcal M_g$ has rank one for $g \geq 3$, as shown by Harer \cite{Ha83} in his study of the homology of the mapping class group. To this day, no purely algebro-geometric proof of this fact is known without Harer’s topological input. Moreover, the rational Picard groups of Hurwitz spaces of branched covers of $\mathbb P^1$ are known only in certain ranges of degree and genus; see \cite{DE96, St00, DP15, Mu23, LL25}.  

Additionally, when $\mu$ is a signature of meromorphic differentials with some simple poles, relation~\eqref{eq:eta-psi} still holds, but the coefficient $m_i+1$ of the corresponding $\psi$-class becomes zero for a simple pole with $m_i = -1$, and hence \eqref{eq:eta-psi} imposes no constraint on that $\psi$-class. The same phenomenon also occurs for strata of $k$-differentials with poles of order $-k$.

Beyond the interior of the strata, one may also study the tautological ring of the compactified moduli space of multi-scale differentials by incorporating the cycle classes of boundary strata; see~\cite[Section~8]{CMZ22}. In genus zero, Devkota \cite{De25} proved that the cohomology ring of the moduli space of multi-scale differentials is generated by the boundary divisors. However, in higher genus, where the boundary structure of the multi-scale space is substantially more intricate than that of $\overline{\mathcal M}_{g,n}$, even conjectural structural descriptions---analogous to Faber’s conjecture and the Faber--Zagier relations for $\overline{\mathcal M}_{g,n}$---remain out of reach for the strata of differentials.

We summarize these discussions in the following question. 

\begin{problem}
    For strata of holomorphic differentials, determine which powers of $\eta$ are nonzero. In particular, for varying strata of holomorphic differentials, determine whether the rational Picard group is of rank one generated by $\eta$. 
    
    For strata of meromorphic differentials with simple poles, determine which powers of the $\psi$-classes associated to simple poles are nonzero. 
    
    Furthermore, describe the tautological ring for strata of $k$-differentials and for the moduli space of multi-scale differentials. 
\end{problem}

\subsection{Affine geometry and complete subvarieties in strata of differentials}
\label{subsec:affine}

Affine varieties in algebraic geometry correspond to Stein manifolds in complex analytic geometry. An affine variety can be viewed as the complement of a hypersurface (i.e., a very ample divisor) in a projective variety. 

Affine varieties possess several notable properties. For example, an affine variety $M$ contains no complete algebraic curves; its homology groups vanish in degrees greater than $\dim_{\mathbb C} M$; and the boundary of any complete algebraic compactification of $M$ is connected. Therefore, it is natural to study the affine geometry of strata of differentials in order to derive corresponding consequences for their topology. 

When $\mu$ is a signature of meromorphic differentials, using the positivity of $\psi$-classes and utilizing the relation 
$$(m_i+1)\psi_i = (m_j+1)\psi_j$$ 
in \eqref{eq:eta-psi} for a zero with $m_i \geq 0$ and a pole with $m_j < 0$, one can produce an ample divisor class on the compactified strata of meromorphic differentials whose support lies entirely on the boundary of the compactification. This implies that the strata of meromorphic differentials are affine varieties, as shown in~\cite{Ch24}, where similar ideas are used to show that the nonvarying strata of holomorphic differentials discussed in Section~\ref{subsec:singularities}, as well as the strata of $k$-differentials of infinite area, are affine varieties. However, it remains unknown whether those varying strata of holomorphic differentials in higher genus are affine varieties. 

As mentioned above, a closely related question is whether strata of differentials can contain complete positive-dimensional subvarieties. Since strata of meromorphic differentials are affine varieties---and hence contain no complete algebraic curves---the remaining case is to determine whether strata of holomorphic differentials may contain complete algebraic curves. 

For the unprojectivized strata $\mathcal H(\mu)$ of holomorphic differentials, Gendron \cite{Ge20} applied the maximum modulus principle to the shortest saddle connections of associated translation surfaces and proved that $\mathcal H(\mu)$ contains no complete curves; an alternative proof using positivity of divisor classes was later given in \cite{Ch23}. However, it remains unknown whether the projectivized strata $\mathbb P\mathcal H(\mu)$ of holomorphic differentials can contain complete algebraic curves. As a special case, the minimal spin components $\mathbb P\mathcal H(2g-2)^{\rm odd/even}$ parametrize (non-hyperelliptic) subcanonical points, i.e., canonical divisors with a unique zero of vanishing order $2g-2$. Harris \cite[p.~413]{Ha84} asked whether there exist complete families in the locus of subcanonical points---a question that has remained open for four decades. 

As a comparison, the moduli space of curves $\mathcal M_g$ contains complete curves for $g \geq 3$, owing to the fact that the complex codimension of the boundary in the Satake compactification of $\mathcal M_g$ is greater than one. It is also known that any complete subvariety of $\mathcal M_g$ has dimension at most $g-2$; see the original proof due to Diaz~\cite{Di84}, as well as later proofs \cite{Lo95, GK09} using different methods. However, a sharp upper bound on the dimensions of complete subvarieties in $\mathcal M_g$ remains unknown; in fact, it is not even known whether $\mathcal M_4$ contains a complete algebraic surface. 

Besides investigating whether a moduli space is affine, one can also study geometrically meaningful stratifications of the moduli space such that each stratum is affine. For the moduli space of curves $\mathcal M_g$, we refer to the survey \cite{Mo14} for an  comprehensive overview of this topic. Here, we briefly mention several results and open questions among many others. Note that $\mathcal M_g$ can be stratified by the loci of $k$-gonal curves for $k = 2, \ldots, g$, known as Arbarello's stratification, where the smallest stratum parametrizes hyperelliptic curves. One may ask whether each $k$-gonal stratum (outside smaller strata) is affine; however, this was shown not to be the case by Arbarello and Mondello \cite{AM12}. More generally, Looijenga asked whether $\mathcal M_g$ admits an affine stratification consisting of $g-1$ layers, and such a stratification is known for $g \leq 5$ due to work of Fontanari and Looijenga \cite{FL08}.  

Returning to the moduli spaces of holomorphic differentials, one can similarly ask about affine stratifications of $\mathbb P\mathcal H(\mu)$. For the minimal strata $\mathbb P\mathcal H(2g-2)$ parametrizing subcanonical points as special Weierstrass points, one variant is to fix the Weierstrass semigroup $H$ and ask whether each locus of subcanonical points $\mathbb P\mathcal H(2g-2)^H$ with the fixed semigroup $H$ is affine. In the case where $H$ is among the most specialized semigroups, i.e., when the corresponding subcanonical points do not further degenerate, it is known that $\mathbb P\mathcal H(2g-2)^H$ is affine \cite[Theorem 8.5]{CY25}; however, the question remains open for general $H$. 

It is worth noting that the locus of Weierstrass points with a fixed semigroup can be identified with the locus of smooth deformations in the miniversal deformation space of the corresponding monomial singularity, as shown by Pinkham \cite[Section 13]{Pi74}. In general, for an isolated complete intersection singularity, the discriminant locus is a hypersurface, and hence its complement parametrizing smooth deformations is affine; see \cite[Theorem 4.8]{Lo84}. However, monomial singularities are often not complete intersections. Indeed, given an arbitrary semigroup, it remains unknown whether the corresponding locus of Weierstrass points on smooth curves is nonempty. 

A related quantity is the Dolbeault cohomological dimension of a variety $M$, defined as the supremum of integers $q$ such that $H^{0,q}_{\overline{\partial}}(M, E)\neq 0$ for some holomorphic vector bundle $E$ on $M$. The Dolbeault cohomological dimension measures how far a variety is from being affine; for example, a smooth affine variety has cohomological dimension zero. We refer to \cite[Appendix A]{Mo14} for an introduction to cohomological dimensions.

For strata of holomorphic differentials $\mathbb P\mathcal H(\mu)$, Mondello \cite{Mo17} constructed an exhaustion function whose prototype is based on the ratio of $\ell_{\rm sys}^2$ and the area of a translation surface, where $\ell_{\rm sys}$ denotes the systole, i.e., the shortest saddle connection or closed geodesic of the translation surface. Using this exhaustion function, Mondello showed that the Dolbeault cohomological dimension of $\mathbb P\mathcal H(\mu)$ is at most $g$; see \cite[Theorem C]{Mo17}. Note that this bound is $g$ away from the sharp cohomological dimension zero for an affine variety, which is essentially due to the fact that the signature of the complex Hessian of the area functional is $(g,g)$; see~\cite[Lemma 3.7(c)]{Mo17}.

We conclude the above discussion with the following question. 

\begin{problem}
    Study complete curves and cohomological dimensions for the projectivized strata of holomorphic differentials $\mathbb P\mathcal H(\mu)$.  Ultimately, determine whether $\mathbb P\mathcal H(\mu)$ is an affine variety. 
\end{problem}

\subsection{Ends of strata of differentials}
\label{subsec:end}

We end our discussion on the topology of strata by addressing  their topological ends. Recall that for a topological space $M$, the number of ends of $M$ is defined as the limit of the number of connected components of $M \setminus K_n$, where ${K_n}$ is an exhaustion of $M$ by compact subsets. Using flat geometry and Rauzy classes, Boissy~\cite{Bo15} showed that every connected component of the strata of holomorphic quadratic differentials (including global squares of one-forms) has exactly one topological end.

Using the multi-scale compactification of strata of differentials (as a smooth complex orbifold), Dozier and Grushevsky~\cite[Lemma 3]{DG25} observed that the connectedness of the boundary of a stratum component implies that it has exactly one end. Using the classification of connected components of strata of Abelian differentials and 
explicit degeneration techniques, they proved that the boundary of any stratum component of holomorphic or meromorphic differentials (of complex dimension at least two) is connected, thereby recovering Boissy’s one-end result in the case of holomorphic one-forms and generalizing it to meromorphic one-forms. In the genus-one case, the connectedness of certain boundary divisors was established by Lee~\cite[Appendix A]{DG25}, building on the work of~\cite{LW25} on connected components of generalized strata. 

Note that the boundary of an affine variety (of dimension at least two) in any complete algebraic compactification is connected~\cite{Go69}. In Section~\ref{subsec:affine}, we have seen that strata of meromorphic differentials are affine, which conceptually explains why every stratum component of meromorphic differentials has connected boundary and hence exactly one end; the same conclusion also holds for strata of $k$-differentials of infinite area; see~\cite{Ch25}.

Nevertheless, since it remains unknown whether strata of holomorphic differentials are affine, this conceptual explanation cannot yet be extended to them. In the case of $k$-differentials of finite area for $k \geq 3$, the lack of a classification of connected components makes it difficult to provide explicit descriptions of boundary connectedness. Furthermore, without the ${\rm GL}_2^+(\mathbb R)$-action in the case $k \geq 3$, even showing that the boundary of every positive-dimensional stratum component of finite-area $k$-differentials is nonempty is a nontrivial problem---one that was only recently resolved by Aygun~\cite[Lemma~3.9]{Ay25} using properties of the area function of flat surfaces. Similar questions may also be posed for linear subvarieties in the strata of holomorphic differentials, where the answers remain largely unknown.  

We conclude this discussion with the following question.  

\begin{problem}
    Without relying on the classification of connected components, provide a conceptual proof that the boundary of every higher-dimensional stratum component of holomorphic differentials is connected. 
    
    Furthermore, investigate the boundary structures of strata of $k$-differentials and of linear subvarieties in the strata of holomorphic differentials, and determine the number of their topological ends. 
\end{problem}

\bibliographystyle{alpha}
\bibliography{biblio}

\end{document}